\documentclass[12pt,reqno]{amsart}
\usepackage{amssymb}

\newtheorem{theorem}{Theorem}[section]
\newtheorem{definition}[theorem]{Definition}
\newtheorem{proposition}[theorem]{Proposition}
\newtheorem{corollary}[theorem]{Corollary}
\newtheorem{lemma}[theorem]{Lemma}

\def \proof {\noindent {\bf Proof.}\ \ }

\def \remark {\noindent {\bf Remark.}\ \ }

\def \R {\mathbb{R}}
\def \Z {\mathbb{Z}}
\def \E {\mathbb{E}}
            
\def \T {\mathbb{T}}
\def \H {\mathcal{H}}     
\def \a {\alpha}

\def \e {\varepsilon}
\def \d {\delta}
\def \D {\Delta}

\def \l {\lambda}
\def \L {\Lambda}
\def \s {\sigma}

\def \< {\langle}
\def \> {\rangle}

\def \pone {{\psi_1}}
\def \ptwo {{\psi_2}}

\def \rank {{\rm rank}}
\def \Rank {{\rm Rank}}

\def \HS {{\rm HS}}
\def \Prob {{\rm Prob}}

\def \id {{\it id}}

\def \dens {{\rm dens}}
\def \supp {{\rm supp}}

\begin{document}
\title {Coordinate Restrictions of Linear Operators in $l_2^n$}
\author {R.Vershynin}
\address{Faculty of Mathematics and Computer Science,
         The Weizmann Institute of Science,
         Rehovot 76100, Israel}
\email{vershyn@wisdom.weizmann.ac.il}
\date{\today}

\begin{abstract}
This paper addresses the problem of improving
properties of a linear operator $u$ in $l_2^n$ 
by restricting it onto coordinate subspaces.
We discuss
how to reduce the norm of $u$ by a random coordinate restriction,
how to approximate $u$ by a random operator with small 
  "coordinate" rank,
how to find coordinate subspaces where $u$ is an isomorphism.
The first problem in this list provides a probabilistic
  extension of a suppression theorem of B.~Kashin and L.~Tzafriri,
the second one is a new look at a result of M.~Rudelson 
  on the random vectors in the isotropic position,
the last one is the recent generalization of the
  Bourgain-Tzafriri's invertibility principle.  
The main point is that all the results are independent of $n$,
the situation is instead controlled by the Hilbert-Schmidt norm
of $u$.
As an application, we provide an almost optimal solution to 
the problem of harmonic density in harmonic analysis, 
and a solution to the reconstruction problem for
communication networks which deliver data with random losses.
\end{abstract}

\maketitle

\section{Introduction}

Linear operators in a finite dimensional Hilbert space $H$
constitute one of the most fundamental classes of operators in 
Functional Analysis.

Recall a classic observation for an arbitrary operator $u$ on $H$.  
There can be found an orthonormal basis in $H$ so that, up to 
an isometry, $u$ is a diagonal operator with
respect to that basis, and its diagonal entries
$s_1  >  s_2  >   \cdots  >  s_N  > 0$
satisfy $\sum s_j^2 = \|u\|_\HS^2$,
where $\|u\|_\HS$ denotes the Hilbert-Schmidt norm of $u$.
This provides enough information on how $u$ acts
with respect to the chosen coordinate structure on $H$.
For instance, restricting $u$ onto an appropriate coordinate
subspace $\R^\s$ we can cut off large s-numbers $s_j$ 
to improve the norm of $u$, or to cut off small s-numbers $s_j$
to nicely approximate $u$ by an operator with a smaller rank.

But what if a coordinate structure on $H$ already exists
and has no relation to $u$? In other words,
can one still improve the properties of a linear operator 
$u$ in $l_2^N$ by restricting it onto coordinate subspaces?

As a first result of this paper, we will compute the norm of 
a random coordinate restriction of $u$. Precisely, we bound
\begin{equation}                           \label{i:expect}
  \E \big\| u \big|_{\R^\s} \big\|
\end{equation}
where $\s$ is a random subset of $\{ 1, \ldots, N \}$
of a fixed cardinality $n \le N$.
The history of the question is the following. 
A known theorem of M.~Talagrand \cite{Ta 95} gives an upper
estimate on \eqref{i:expect} for a linear operator from $l_2^N$ 
into a Banach space $X$ (see \cite{Ta 98} for a further generalization).
However, the estimate of M.~Talagrand is close to being sharp
only in a certain, quite restrictive, range of $n$
(needed in applications to the $\L(p)$-problem).
After partial results of B.~Kashin \cite{Ka} and A.~Lunin \cite{Lu}, 
B.~Kashin and L.~Tzafriri \cite{Ka-Tz} 
produced an argument which essentially proves an optimal 
estimate on the {\em minimum} of $\big\| u \big|_{\R^\s} \big\|$
over all subsets $\s$ as above. It equals to 
\begin{equation}                                \label{i:v}
  \min_{|\s| = n} \big\| u \big|_{\R^\s} \big\|
    \le  C \Big( \sqrt{ \frac{n}{N} } + \sqrt{ \frac{h}{N} } \Big)
\end{equation}
where $h = \|u\|_\HS^2$ and $u$ is assumed norm one (see \cite{V}). 
The crucial step in the proof of \eqref{i:v} is Grothendieck's
factorization (see \cite{Le-Ta} Proposition 15.11),
which gives no information about $\s$ except that it can be 
found in a random subset of $\{1, \ldots, N \}$ 
of cardinal, say, $2 h$. 
Note that \eqref{i:v} is equivalent to 
\begin{equation}                                     \label{minmax}
\min_{|\s| = n} \big\| u \big|_{\R^\s} \big\|
  \le  C \Big( \sqrt{ \frac{n}{N} } + \max_j \|u e_j\| \Big),
\end{equation}
because by Chebyshev's inequality at least $\frac{1}{2} N$ 
numbers $\|u e_j\|$ do not exceed $2 \sqrt{ \frac{h}{N} }$.
We will prove that, up to a logarithmic factor, 
the same estimate holds for a random set $\s$.

\begin{proposition}                            \label{propsup}
  Let $u$ be a norm one linear operator in $l_2^N$.
  Consider an integer $1 \le n \le N$. Then for random subset 
  $\s$ of $\{ 1, \ldots, N \}$ of cardinal $|\s| = n$ 
  $$
  \E \big\| u \big|_{\R^\s} \big\|
  \le  C \log{n} 
      \Big( \sqrt{ \frac{n}{N} } + \max_j \|u e_j\| \Big).
  $$
\end{proposition}

In particular, for the threshold dimension, $|\s| = h$, 
we have
\begin{equation}                                      \label{i:sup}
\E \big\| u \big|_{\R^\s} \big\|
  \le  C \log{h} \cdot \max_j \|u e_j\|.
\end{equation}
Note that the dimension $N$ plays no role in \eqref{i:sup}. 
Instead, the situation is completely controlled by the
parameter $h = \|u\|_\HS^2$.

The key to the proof of Proposition \ref{propsup} is 
the non-commutative Khinchine inequality in the Schatten 
class $C_p^n$, due to F.~Lust-Piquard and G.~Pisier
(see \cite{P}). Its usefulness for coordinate restrictions
in $l_2^N$ was recognized by G.~Pisier. This provided 
an alternative approach to a lemma of 
M.~Rudelson \cite{R} (see also \cite{P})
previously proved by a delicate construction
of a majorizing measure. 
We will use the following non-simmetric version of 
M.~Rudelson's lemma, which also follows from the 
non-commutative Khinchine inequality.
For a finite set of vectors $x_j, y_j$ in $\R^n$
\begin{align}                    \label{i:xjyj}
\E \Big\| \sum_j  \e_j x_j \otimes y_j \Big\|
  \le C \sqrt{\log n}  
      \bigg( 
        &\max_j \|x_j\|
           \cdot \Big\| \sum_j  y_j \otimes y_j \Big\|^{1/2}  \notag \\
      + &\max_j \|y_j\|
           \cdot \Big\| \sum_j  x_j \otimes x_j \Big\|^{1/2}
      \bigg)
  \end{align}
where $\e_j$ are independent Bernoulli random variables with 
$\Prob \{ \e_j = 1 \}  =  \Prob \{ \e_j = -1 \}  =  \frac{1}{2}$.

Next, we will prove an "approximation" counterpart
of Proposition \ref{propsup}.
Recall a well known inequality for the approximation numbers of 
an operator $u$ in $l_2^N$:
\begin{equation}                                  \label{i:an}
  a_n(u)   \le   \sqrt{ \frac{h}{n} }
\end{equation}    
where $a_n  =  \inf \big\{ \|u - u_1\| : \; \rank u_1 < n \big\}$.
We will obtain a coordinate version of \eqref{i:an}
by a different look on arguments of M.~Rudelson \cite{R}.

\begin{theorem}                \label{i:marks}
  Let $u$ be a norm one linear operator in $l_2^N$, 
  and $h  =  \|u\|_{HS}^2$.
  Then for any integer $n > 1$ there exists a diagonal operator 
  $\D$ in $l_2$ such that $\rank \D \le n$ and
  \begin{equation}                                 \label{i:appr}       
    \| u (\D - \id) u^* \|  
    \le  C \sqrt{\log n} \cdot \sqrt{ \frac{h}{n} }.
  \end{equation}
\end{theorem}

Examples show that both $u$ and $u^*$ are needed in \eqref{i:appr}.
$\D$ is a random diagonal operator whose entries are multiples
of independent selectors. It depends only on the values of $\|u e_j\|$;
the larger $\|u e_j\|$ is, the more likely the $j$-th entry of $\D$
is not zero.
For such random operator $\D$,  \eqref{i:appr} holds 
with large probability. Namely, if we set 
$\e = C \sqrt{\log n} \cdot \sqrt{ \frac{h}{n} }$
then 
\begin{equation}                     \label{i:conc}
  \Prob \big\{ \| u (\D - \id) u^* \| > t \e  \big\}   
      \le  3 \exp (-t^2).
\end{equation}

The reader should note that the dimension $N$ of the space 
plays no role in this result as well as in Proposition \ref{propsup},
and the situation is again controlled by the parameter 
$h = \|u\|_\HS^2$.
This phenomenon seems quite general. 
It roughly says that for a linear operator $u$ in $l_2^N$
the Hilbert-Schmidt norm of $u$ (an not the rank, for example)
is responsible for the essential properties of $u$.
Another instance of this phenomenon is the recent extension
of Bourgain-Tzafriri's principle of restricted invertibility
\cite{V}.

\begin{theorem}                                   \label{i:rip}
  Let $u$ be a norm one linear operator in $l_2^N$, 
  and $h  =  \|u\|_{HS}^2$.
  Then for any $\e > 0$ there exists a subset $\s$ 
  of $\{ 1, \ldots, N \}$ of cardinal $|\s| > (1 - \e) h$
  so that the sequence $(T e_j)_{j \in \s}$ is $C(\e)$-equivalent 
  to an orthogonal basis.
\end{theorem}

In other words, under the hypotheses of Theorem \ref{i:rip}
there exists an isomorphism $T$ with $\|T\| \|T^{-1}\| < C(\e)$
which takes $e_j$ to $u e_j / \|u e_j\|$ for $j \in \s$.
We see again that the parameter $h = \|u\|_\HS^2$ governs
the situation, and the dimension $N$ is unimportant.

There is an important application of Theorem \ref{i:rip}
to "unbounded" operators $u$ in $l_2^N$. It says that 
if the norms $\|u e_j\|$ are controlled, 
then $u$ is a nice isomorphism on a large coordinate subspace. 

\begin{corollary}                                 \label{i:cor}
  Let $u$ be a linear operator in $l_2^N$ such that 
  $\|u e_j\| = 1$ for all $j = 1, \ldots, N$.
  Then for any $\e > 0$ there exists a subset $\s$ 
  of $\{ 1, \ldots, N \}$ of cardinal 
  $|\s| > (1 - \e) \frac{N}{\|u\|^2}$
  so that
  \begin{equation}                          \label{i:iso}
    c_1(\e) \|x\|  \le  \|u x\|  \le  c_2(\e) \|x\|
    \ \ \ \ \text{for $x \in \R^\s$}.
  \end{equation}
\end{corollary}

\noindent This is a direct generalization of a theorem 
of J.~Bourgain and L.~Tzafriri, who proved Corollary \ref{i:cor} for 
some $0 < \e < 1$ and with only the lower bound in 
\eqref{i:iso}.
The upper bound is also nontrivial, as we do not assume 
the operator $u$ to be well bounded. 

An application of Corollary \ref{i:cor} to harmonic analysis
generalizes results on the problem of harmonic density
\cite{B-Tz}.
Let $T$ be the circle with the normalized Lebesgue measure $\nu$,
and $B$ is a subset of $T$ of positive measure.
The two norms naturally arise here,
$$
\| f \|_{L_2(B)}  
  =  \Big( \frac{1}{\nu(B)} 
           \int_B |f|^2 \; d\nu
     \Big)^{1/2}
\ \ \ \text{and} \ \ \ 
\| f \|_{L_2(T)}  
  =  \Big( \int_\T |f|^2 \; d\nu
     \Big)^{1/2}.
$$ 
In general there is no relation between $\|f\|_{L_2(B)}$ 
and $\|f\|_{L_2(T)}$. However, suppose $B$ is a half-circle; 
then it is easily seen that the two norms are equal for the functions
$f$ whose Fourier transform is supported by $2 \Z$.
Then (a rather vague) question is 
$$
\begin{aligned}
  \qquad&\text{For what functions $f$ are the two norms} \\
        &\text{$\|f\|_{L_2(B)}$ and $\|f\|_{L_2(\T)}$ equivalent?}
\end{aligned}
$$
More specifically, consider functions $f$ whose 
Fourier transform $\hat{f}$
is supported by a fixed set of integers $\L$.
How dense can $\L$ be so that the two norms are still
equivalent?

This problem in a weaker form was stated by W.~Schachermayer.
His question was whether there exists a set $\L$ of integers 
such that $f$ does not vanish a.e. on $B$ provided 
$\supp \hat{f}  \subset \L$.
Answering this question positively, J.~Bourgain and L.~Tzafriri
proved that the existence of $\L$ with density
$c \nu(B)$ for which 
\begin{equation}                             \label{i:bigger}
  \|f\|_{L_2(B)}  \ge  c \|f\|_{L_2(\T)}
\end{equation}
whenever $\supp \hat{f}  \subset \L$.
The proof relies on the principle of restricted invertibility.
The extension of this principle, Corollary \ref{i:cor}, 
can be used to prove the
reverse inequality in \eqref{i:bigger} and also to obtain a 
nearly optimal density $(1 - \e) \nu(B)$ of $\L$. 

\begin{theorem}                                  \label{harmonic}               
  Let $B$ be a subset of $\T$ of positive Lebesgue 
  measure, and $\e > 0$. 
  Then there exists a set of integers $\L$ with 
  (two-sided) density $\dens \L  >  (1 - \e) \nu(B)$ 
  so that
  \begin{equation}                               \label{equivalence}
    c_1(\e) \|f\|_{L_2(\T)} 
    \le \|f\|_{L_2(B)} 
    \le c_2(\e) \|f\|_{L_2(\T)}
  \end{equation}
  whenever the Fourier transform of $f$ is supported by $\L$.
\end{theorem}  

The two-sided density here is
$\dens \L = \lim_{n \to \infty} \frac{ |\L \cap [-n, n]| }{2 n}$.
This definition seems to be more natural for subsets of $\Z$
than the usual "one-sided" definition. 

In Appendix, we return to Theorem \ref{i:marks} and 
discuss its rather unexpected application. 
Together with the concentration inequality \eqref{i:conc}, 
Theorem \ref{i:marks} provides a solution to the reconstruction 
problem in communication systems such as the Internet
which deliver data with random losses \cite{G-K}. 

\thanks I am grateful to P.~Casazza and to G.~Schechtman for useful 
discussions.

\section{Suppressions on Coordinate Subspaces}           \label{suppressions}

\subsection{Result}                                  

In this section we will compute the norm of the random coordinate 
restriction $\big\| u \big|_{\R^\s} \big\|$
of an arbitrary linear operator $u$ on $l_2^N$.
The size of $\s$ is a fixed integer $n$
not necessarily equal $h ( = \|u\|_\HS^2 )$ as 
in Proposition \ref{propsup}.

The {\em minimum} of $\big\| u \big|_{\R^\s} \big\|$
over all subsets $\s$ of cardinal $n$ is provided by 
Kashin-Tzafriri's theorem \cite{Ka-Tz} (see \cite{V} for a proof).

\begin{theorem}                                   \label{ktztheorem}
  (B.~Kashin, L.~Tzafriri). \
  Let $u$ be a norm one linear operator in $l_2^N$,
  and $h = \|u\|_\HS^2$.
  Then for any integer $n > 1$ 
  there exists a subset $\s$ of $\{ 1, \ldots, N \}$
  of cardinal $|\s| = n$ such that
  \begin{equation}                        \label{urs}
    \big\| u \big|_{\R^\s} \big\|
    \le  c \Big( \sqrt{ \frac{n}{N} } + \sqrt{ \frac{h}{N} } \Big).
  \end{equation}
\end{theorem}

We see that the treshold value for the size of $\s$  
is $h$:
for $n \le  h$ the best restriction 
$\big\| u \big|_{\R^\s} \big\|$
is bounded by $\sqrt{ \frac{h}{N} }$, while 
for a larger size $n$ the best bound is $\sqrt{ \frac{n}{N} }$.
This is illustrated by an example in Section \ref{form},
which yields that \eqref{urs} is sharp up to a constant.

Before we state a probabilistic extension of Theorem \ref{ktztheorem},
let us specify what we mean by a random subset $\s$. 
There are several equivalent definitions, of which 
the following one seems more convenient throughout
the present paper.
Let $A$ be a finite set and $0 < \d < 1$.
Consider a subset $\s$ of $A$ by taking (rejecting) 
each element of $A$ independently with probability $\d$
(respectively, $1 - \d$). Then the cardinal of $\s$ is
concentrated around $n = \d |A|$.
We then call $\s$ a {\em random subset} 
of $A$ of cardinal $|\s|  \sim  n$.

\begin{theorem}                         \label{suppression}
  Let $u$ be a norm one linear operator in $l_2^N$.
  Let $M = \max_j \|u e_j\|$.
  Then for any integer $n > 1$ 
  and for a random subset $\s$ of $\{ 1, \ldots, N \}$ 
  of cardinal $|\s|  \sim  n$
  \begin{equation}                      \label{eqsup}
    \E \|u \big|_{\R^\s} \|  
      \le  C \sqrt{\log n} \Big( \sqrt{ \frac{n}{N} } 
                                    + \sqrt{\log n} \; M 
                           \Big).
  \end{equation}
\end{theorem}

The key to the proof is a non-simmetric version of 
M.~Rudelson's lemma from \cite{R} (see also \cite{P}).

\begin{lemma}                              \label{Pisier}
  Let $x_j, y_j$ be a finite set of vectors in $\R^m$. 
  Then 
  \begin{align}                    \label{xjyj}
  \E \Big\| \sum_j  \e_j x_j \otimes y_j \Big\|
    \le C \sqrt{\log m} 
         \bigg( &\max_j \|x_j\|
            \cdot \Big\| \sum_j  y_j \otimes y_j \Big\|^{1/2}  \notag \\
          + &\max_j \|y_j\|
            \cdot \Big\| \sum_j  x_j \otimes x_j \Big\|^{1/2}
         \bigg) .
  \end{align}
\end{lemma}

\noindent This lemma is a consequence of the non-commutative 
Khinchine inequality in the Schatten space $C_p^m$, with
optimal constant $O(\sqrt{p})$, for $p = \log m$.
This inequality is a result of F.~Lust-Piquard and G.~Pisier
(see \cite{P}).

\begin{theorem} ( F.~Lust-Piquard,  G.~Pisier). \ \        \label{noncom}
  Assume $2 \le p < \infty$.
  Then there is a constant $B_p  \le  C \sqrt{p}$ such that
  for any finite sequence $(X_j)$ in $C_p^m$
  $$
  R(X_j) 
    \le \Big\| \sum \e_j X_j \Big\|_{L_p(C_p^m)}    
    \le  B_p \cdot R(X_j),
  $$
  where 
  $$
  R(X_j)  =  \Big\| \Big( \sum X_j^* X_j \Big)^{1/2} \Big\|_{C_p^m}
   \vee \Big\| \Big( \sum X_j X_j^* \Big)^{1/2} \Big\|_{C_p^m}.
  $$
\end{theorem}

\noindent {\bf Proof of Lemma \ref{Pisier}.\ }
Note that for $p = \log m$
\begin{equation}                            \label{cpn}
  \|X\|  \le  \|X\|_{C_p^m}  \le  e \|X\|.
\end{equation}
Then we can apply Theorem \ref{noncom} for $X_j  =  x_j \otimes y_j$
noting that 
$X_j^* X_j  =  \|y_j\|^2 x_j \otimes x_j$ and 
$X_j X_j^*  =  \|x_j\|^2 y_j \otimes y_j$.
We get 
\begin{align*}
\E \Big\| \sum_j  \e_j x_j \otimes y_j \Big\|
  &\le  \E \Big\| \sum_j  \e_j x_j \otimes y_j \Big\|_{C_p^m}  \\
  &\le  \Big\| \sum_j  \e_j x_j \otimes y_j \Big\|_{L_p(C_p^m)}  \\
  &\le  C \sqrt{p} \cdot 
        \bigg( \Big\| \Big( \sum  \|y_j\|^2 x_j \otimes x_j 
                      \Big)^{1/2}
               \Big\|_{C_p^m}     \\
    &\qquad  + \Big\| \Big( \sum  \|x_j\|^2 y_j \otimes y_j 
                      \Big)^{1/2}
               \Big\|_{C_p^m}    
         \bigg).
\end{align*}
By \eqref{cpn} we can replace both $\| \cdot \|_{C_p^m}$-norms
by $\| \cdot \|$-norms, which easily leads to the completion 
of the proof.
\endproof

\noindent {\bf Proof of Theorem \ref{suppression}.\ }
For $x_j = u e_j$, $j = 1, \ldots, N$ we can write
$$
u  =  \sum_{j \le N} e_j \otimes x_j.
$$
We need to compute
$$
E := \E \| u \big|_{\R^\s} \|
   = \E \Big\| \sum_{j \le N} \d_j e_j \otimes x_j \Big\|,
$$
where $\d_j$ are $\{ 0, 1 \}$-valued independent random variables
with $\E \d_j = \d = \sqrt{ \frac{n}{N} }$.
This will be done by a usual symmetrization and applying
\eqref{xjyj} twice; first to bound $E$ in terms of 
$E_1  =  \E \Big\| \sum \d_j x_j \otimes x_j \Big\|$,
and then again to compute $E_1$.

Now we pass to a detailed proof. 
The standard symmetrization procedure (see \cite{Le-Ta} Lemma 6.3)
yields
\begin{align*}
E &\le  \E \Big\| \sum_{j \le N} 
                  (\d_j - \d) e_j \otimes x_j 
           \Big\|
        +  \d \|u\|          \\
  &\le  2 \E \E_\e 
          \Big\| \sum_{j \le N} \e_j \d_j e_j \otimes x_j 
          \Big\|
        +  \d.
\end{align*} 
Then we apply \eqref{xjyj} to bound 
$\E_\e \Big\| \sum \e_j  e_j \otimes ( \d_j x_j ) \Big\|$.
Clearly in \eqref{xjyj} we can set $m$ equal 
$$
N(\d) :=  e \vee \Big( \sum_{j \le N} \d_j \Big).
$$
Then using Cauchy-Schwartz and Jensen inequalities we obtain
\begin{align}                     \label{eless}
E &\le  C \Big( \E \log N(\d) \Big)^{1/2}
        \left[
          \E \Big( \Big\| \sum_{j \le N} \d_j x_j \otimes x_j 
                   \Big\|^{1/2} + M
             \Big)^2
        \right]^{1/2}
        + \d                    \notag \\
  &\le  C \Big( \log \E N(\d) \Big)^{1/2} 
        \left[
          \Big( \E \Big\| \sum_{j \le N} \d_j x_j \otimes x_j \Big\|
          \Big)^{1/2} + M
        \right]
        + \d                    \notag \\
  &\le  C \sqrt{\log \d N}
        \big( E_1^{1/2} + M \big) 
        + \d, 
\end{align}     
where 
$$
E_1  =  \E \Big\| \sum_{j \le N} \d_j x_j \otimes x_j 
           \Big\|.
$$
Therefore the problem reduced to computing $E_1$.
By the standard symmetrization and noting that
$ \big\| \sum x_j \otimes x_j \big\|
  =  \| u u^* \|  \le  1$
we have
\begin{align*}
E_1 
  &\le \E \Big\| \sum_{j \le N} 
                  (\d_j - \d) x_j \otimes x_j 
          \Big\|
        +  \d          \\
  &\le  2 \E \E_\e 
          \Big\| \sum_{j \le N} \e_j \d_j x_j \otimes x_j 
          \Big\|
        +  \d. 
\end{align*}
We apply \eqref{xjyj} again.
\begin{align*}                 
E_1
  &\le  C \Big( \E \log N(\d) \Big)^{1/2}
        \cdot M 
        \cdot \Big( \E \Big\| \sum_{j \le N} \d_j x_j \otimes x_j 
                       \Big\|
             \Big)^{1/2}
        + \d                   \\
  &\le  C \sqrt{\log \d N}
        \cdot M
        \cdot E_1^{1/2}  
        + \d. 
\end{align*} 
It follows that 
$$
E_1^{1/2}  \le  C \big( \sqrt{\log \d N}
                        \cdot M 
                        + \sqrt{\d}
                  \big).
$$
Combining this estimate with \eqref{eless} we obtain
$$
E \le  C \sqrt{\log \d N}
       \big( \sqrt{\log \d N} \cdot M  + \sqrt{\d}
       \big).
$$
This completes the proof.
\endproof

\subsection{Dimension Independent Corollary}                  

Let us take a closer look at \eqref{eqsup}.
If we disregard for a moment the distracting logarithmic terms
we will see the point of this estimate.
It essentially says that the norm of $u$ is bounded
on most coordinate subspaces of dimension $n$ by 
the maximum of $\sqrt{ \frac{n}{N} }$ and $M  =  \max_j \|u e_j\|$.
The size of $M$ is a natural obstacle here;
the only a priori information about $M$ is that
\begin{equation}                               \label{M}
  \sqrt{ \frac{h}{N} }  \le  M  \le  1
\end{equation}
where, as usual, $h = \|u\|_\HS$.
However, the upper bound for $M$ can be reduced 
to essentially $\sqrt{ \frac{h}{N} }$ as explained in the 
introduction. 

So, the new information \eqref{eqsup} gives is that 
the optimal estimate of Kashin and Tzafriri on the minimal
suppression $\big\| u \big|_{\R^\s} \big\|$ 
essentially holds for {\em most} subsets $\s$. 

Taking $n = h \log h$, we obtain through \eqref{M}
$$                                  
\E \big\| u \big|_{\R^\s} \big\|
   \le  C \sqrt{\log h} 
        \Big( \sqrt{ \frac{h \log h}{N} }
            + \sqrt{\log h} \; M
        \Big)     
    =   C \log h \cdot M,
$$
and this estimate seems to be sharp up to a constant. 
We single it out as a corollary.

\begin{corollary}                                        \label{loglog}     
  Let $u$ be a norm one linear operator in $l_2^N$, 
  and $h  =  \|u\|_{HS}^2$.  
  Then for a random subset $\s$ of $\{ 1, \ldots, N \}$
  of cardinal $|\s|  \sim h  \log h$ 
  \begin{equation}                               \label{newlogh}
    \E \big\| u \big|_{\R^\s} \big\|
    \le  C \log h \cdot \max_j \|u e_j\|.
  \end{equation}
\end{corollary}

Quite remarkably, the dimension $N$ is not important here. 
We obtained a probabilistic version of Kashin-Tzafriri's estimate  
\begin{equation}                             \label{i:ktz}
  \min_{ |\s| = h } \big\| u \big|_{\R^\s} \big\|
  \le  c \sqrt{ \frac{h}{N} }
\end{equation}
since $M = \max_j \|u e_j\|$ 
can be easily reduced to $c \sqrt{ \frac{h}{N} }$ as explained above.  

How sharp is Theorem \ref{suppression}?
It can be shown (considering the example from Section \ref{form})
that the first logarithmic term in \eqref{eqsup} is necessary, 
and therefore \eqref{eqsup} is sharp up to the second logarithmic
factor (before $M$).
This shows in particular that Theorem \ref{suppression} is sharp 
for large $n$.

\section{Approximation on Coordinate Subspaces}      \label{global}

\subsection{Result}

In the present section we will prove Theorem \ref{i:marks}
and the $\ptwo$-estimate \eqref{i:conc}.
To highlight the diension independentness in Theorem \ref{i:marks},
we state it as an infinite dimensional result.

\begin{theorem}                \label{marks}
  Let $u$ be a norm one linear operator in $l_2$,
  and $h  =  \|u\|_{HS}^2  <  \infty$.
  Then for any integer $n > 1$ there exists a diagonal operator 
  $\D$ in $l_2$ such that $\rank \D \le n$ and
  \begin{equation*}                    
    \| u (\D - \id) u^* \|  
    \le  C \sqrt{\log n} \cdot \sqrt{ \frac{h}{n} }.
  \end{equation*}
\end{theorem}

Setting $C(\e) = C \e^{-2} \log(1/\e)$, we 
obtain an "approximation" counterpart of Corollary \ref{loglog}.

\begin{corollary}       
  Let $u$ be a linear operator in $l_2$ with
  $\|u\| = 1$  and $h = \|u\|_{HS}^2 < \infty$.
  Then for any $\e > 0$ there exists a diagonal operator 
  $\D$ in $l_2$ 
  such that
  \begin{equation}                     \label{lessthane}
    \| u (\D - \id) u^* \|  < \e,
  \end{equation}
  and $\rank \D  \le  C(\e) h (1 + \log h)$.
\end{corollary}

This result can be viewed as a coordinate version of 
the basic inequality for the approximation numbers of $u$
which follows from \eqref{i:an},
$$                                
a_{c(\e) h}(u)  \le  \e
$$
where $c(\e) = \e^{-2}$.

\qquad

The diagonal operator $\D$ is random, and its diagonal entries 
$\D(j)$ can be easily described.
Let $K = n / h$.
\begin{itemize}
  \item {If $\|u e_j\| = 0$ then we set $\D(j) = 0$.}
  \item {If $K \|u e_j\|^2  > 1$ then we set $\D(j) = 1$.}
  \item {If $0 < K \|u e_j\|^2  \le  1$ 
         then $\D(j)$ is a random variable independent 
         of the other entries and distributed as
         \begin{eqnarray}             \label{distribution}
           \Prob \Big\{ \D(j) = \frac{1}{K \|u e_j\|^2} \Big\}  
           =  1 - \Prob \big\{ \D(j) = 0 \big\}
           =  K \|u e_j\|^2.
         \end{eqnarray}
        }
\end{itemize}

\qquad

\noindent {\bf Proof of Theorem \ref{marks}.\ }
The proof requires just a different look at the argument of 
M.~Rudelson \cite{R}. The key step is the application 
of Lemma \ref{Pisier} for $x_j = y_j = u e_j$.
 
We are going to prove for the random operator $\D$ that 
\begin{eqnarray}                           \label{e}
  E := \E \| u (\D - \id) u^* \|  \le  \e.
\end{eqnarray}
First note that $\D$ has finite rank with probability $1$, 
because by the Chebyshev inequality
\begin{eqnarray}                            \label{rankd}
\E \rank \D
  &\le&  | \{ j : K \|u e_j\|^2  \ge  1 \} |
         +  \sum_{j \ge 1} K \|u e_j\|^2  \notag\\
  &\le&  K h  +  K h
     <   \infty.
\end{eqnarray} 
This observation allows us to concentrate only on 
finite-dimensional operators $u$.
To make this claim precise, put $u_N = u P_N$,
where $P_N$ is the coordinate projection in $l_2$
onto $\R^N$.
Let $\D_N$ be the diagonal operator defined as above for the 
operator $u_N$. Then $\D_N = \D P_N = \D$
with probability$\to 1$ as $N \to \infty$. 
Thus
\begin{eqnarray*}
u (\D - \id) u^*  -  u_N (\D_N - \id) u_N^* 
  &=&  u (\D - \id) u^*  - {} \\
       & & {} - u P_N (\D P_N - \id) P_N u^* \\ 
  &=&  (u \D u^* - u \D P_N u^*) + {} \\
       & & {} + (u u^* - u P_N u^*)  \\
  &=&  u u^* - u P_N u^*
\end{eqnarray*}
with probability $\to 1$ as $N \to \infty$.
The norm of this operator vanishes as $N \to \infty$ because
$\|u\|_\HS < \infty$.
This shows that the difference between
$\|u (\D - \id) u^*\|$ and $\|u_N (\D_N - \id) u_N^*\|$
is at most $\e_N$, with probability at least
$1 - \e_N$, where $\e_N \to 0$ as $N \to \infty$. 
Therefore we can assume in (\ref{e}) that $u$ acts in 
a finite dimensional space $l_2^N$.
 
As we already saw in (\ref{rankd}), 
the operator $\D$ defined has the required rank
\begin{equation}                              \label{erankd}
  \E \rank \D  \le  2 K h  =  2 n
\end{equation} 
(the factor $2$ is of course not important).
Let
$$
x_j = u e_j, \ \ j = 1, \ldots, N.
$$
Then we can write 
$$
u u^*  =  \sum_{j=1}^N  x_j \otimes x_j,
$$
so
\begin{eqnarray}                       \label{expect}
E &=&  \E \Big\| \Big( 
                   \sum_{j \le N} \D(j) x_j \otimes x_j 
                 \Big) 
                 - u u^* 
          \Big\|          \notag\\
  &=&  \E \Big\| \sum_{j \le N} (\D(j) - 1) \ x_j \otimes x_j 
          \Big\|.
\end{eqnarray}
At this point we can assume that 
$$
0  <  K \|x_j\|^2  <  1
\ \ \ \ \text{for all $j = 1, \ldots, n$.}
$$
Indeed, if $K \|x_j\|^2$ is either $0$ or $1$, then 
by the construction the summand
$(\D(j) - 1) \ x_j \otimes x_j$
vanishes and contributes nothing to the sum in (\ref{expect}).
Therefore, we can assume that in (\ref{expect}) all
$\D(j)$'s are independent random variables with distribution
(\ref{distribution}).
Since $\E (\D(j) - 1) = 0$ for each $j$, 
we can apply the standard symmetrization procedure 
(see \cite{Le-Ta} Lemma 6.3) which gives
$$
E  \le  2 \E \E_\e \Big\| \sum_{j \le N}  
   \e_j \D(j) x_j \otimes x_j  \Big\|,
$$
where the expectation $\E_\e$ is taken according to the 
Rademacher variables $\e_j$.

To bound the latter expectation, we fix a realization of $\D$ 
and apply Lemma \ref{Pisier} 
for $x_j = y_j = \D(j)^{1/2} x_j$, $j = 1, \ldots, N$. 
The number of non-zero terms in this sequence is 
at most $\rank \D$, so we can assume 
that $m$ in the lemma equals 
$$
\Rank \D  :=  e \vee \rank \D.
$$
We obtain
$$
E  \le  C \E \left[ ( \log \Rank \D )^{1/2}
              \cdot \Big( \max_{j \le N} \D(j)^{1/2} \|x_j\| \Big)
              \cdot \Big\| \sum_{j \le N} \D(j) x_j \otimes x_j 
                    \Big\|^{1/2}
              \right].     
$$
Note that 
$$
\max_{j \le N} \D(j)^{1/2} \|x_j\| = K^{-1/2}.
$$
Then by the Cauchy-Schwartz inequality, Jensen's inequality
and using 
$\big\| \sum_{j=1}^N x_j \otimes x_j \big\|  
  =  \|u u^*\|  \le  1$
we obtain through (\ref{erankd})
\begin{align}                                                 \label{E}
E &\le  C K^{-1/2} 
        \Big( \E \log \Rank \D \Big)^{1/2} 
        \Big( \E \Big\| \sum_{j=1}^N \D(j) x_j \otimes x_j 
                    \Big\| \Big)^{1/2}     \notag\\
  &\le  C K^{-1/2} 
        \Big( \log \E \Rank \D \Big)^{1/2} 
        (E + 1)^{1/2}                      \notag\\
  &\le  C K^{-1/2} 
        (\log n)^{1/2} 
        (E + 1)^{1/2}                      \notag\\
  &=    C \sqrt{\log n} \sqrt{ \frac{h}{n} }
        (E + 1)^{1/2}
\end{align}                       
We can assume that 
$C \sqrt{\log n} \sqrt{ \frac{h}{n} }  \le  1$,
otherwise the conclusions of the theorem hold 
simply with $\D = 0$.
Then \eqref{E} implies
$$
\E  \le  C \sqrt{\log n} \sqrt{ \frac{h}{n} }.
$$
Hence the conclusions of the theorem hold with
non-zero probability. 
This completes the proof.
\endproof

\subsection{Remarks about the form of approximation}          \label{form}

Our first observation is that both $u$ and $u^*$ 
are needed in \eqref{lessthane},
whatever the norm $\|u\|_\HS$ is.
To see this, consider the following example 
essentially borrowed from \cite{Ta 95}.
We consider two positive integers $h$ and $k$ and set $N = h k$.
We define an operator $u$ in $l_2^N$ "blockwise" by its 
action on the coordinate basis:
$$
(u e_i)_{j \le N} := 
\Big(
\underbrace{
  \overbrace{ \frac{e_1}{\sqrt{k}}, \ldots, \frac{e_1}{\sqrt{k}}
            }^{\text{$k$ vectors}},
  \overbrace{ \frac{e_2}{\sqrt{k}}, \ldots, \frac{e_2}{\sqrt{k}}
            }^{\text{$k$ vectors}},
  \ldots, 
  \overbrace{ \frac{e_h}{\sqrt{k}}, \ldots, \frac{e_h}{\sqrt{k}}
            }^{\text{$k$ vectors}}
}_{\text{$N$ vectors}}
\Big).
$$
In other words, we divide $\{1, \ldots, N \}$ into $h$ 
sets $A_l$, $l \le h$, of cardinal $k$, 
and for $j \in A_l$ we set 
$u e_j  =  \e_l / \sqrt{k}$.
Then obviously $\|u\| = 1$ and $\|u\|_{\HS}^2 = h$.
Let $\D$ be an arbitrary diagonal operator in $l_2^N$ with
$\rank \D = n$, where $n = C(\e) h \log h$.
Consider $\s  =  \{ i, \ 1 \le i \le N, \ \D(i) = 0 \}$. 
Since $|\s| \ge N - n$, there exists an integer $l \le h$
such that
$$
| \s \cap A_l |  \ge  \frac{N - n}{h}.
$$
Notice that for $j \in \s \cap A_l$ all vectors 
$u (\id - \D) e_j$ equal the same vector $e_l / \sqrt{k}$.
This implies that 
$$
\| u (\id - \D) \|
  \ge \sqrt{ \frac{ | \s \cap A_l | }{k} }
  \ge \sqrt{ \frac{N - n}{h k} }
   =  \sqrt{ \frac{N - n}{N} }
$$
If $k$ is chosen large enough, then the ratio 
$\frac{n}{N} = \frac{C(\e) \log h}{k}$ is small.
In particular
$$
\| (\id - \D) u^* \| = \| u (\id - \D) \|  \ge  1/2.
$$
This contradicts \eqref{lessthane}.

This example can also be used to show that Theorem \ref{ktztheorem}
is sharp up to a constant, and that Theorem \ref{suppression}
is sharp up to the second logarithmic factor
(we leave this to the interested reader).

\qquad

A comment is in order about the form of the operator $\D$.
One is tempted to say that $\D$ should look like a
projection, i.e. $\D  =  \a P$
for some number $\a$ and a coordinate projection $P$.
This is not always true (again whatever the norm $\|u\|_\HS$ is),
as the following modification of the previous example shows.
We again consider two positive integers $h$ and $k$,
but set $N - 1 = h k$.
We define an operator $u$ in $l_2^N$ similarly to 
the previous example. 
Namely, we divide $\{1, \ldots, N - 1 \}$ into $h$ 
sets $A_l$, $l \le h$, of cardinal $k$, 
and for $j \in A_l$ we set 
$$
x_j  =  u e_j  =  \e_l / \sqrt{k}.
$$
Finally, set 
$$
u e_N = e_{h+1}.
$$
Note that 
\begin{eqnarray}                       \label{ustaru}
u u^*  
  &=&  \Big( \sum_{i \le N-1} x_i \otimes x_i \Big)  
         +  e_{n+1} \otimes e_{h+1}   \notag \\
  &=&  \sum_{j \le h+1} e_j \otimes e_j.
\end{eqnarray}
Then $\|u\| = 1$ and $\|u\|_{HS}^2 = h+1$. 
Let $n = C(\e) h \log h$.

Now assume that there exists a number $\a$ and a set 
$\s \subset \{ 1, 2, \ldots, N \}$, $|\s| \le n$,  so that 
\begin{equation}                        \label{difference}
  \| u (\a P_\s - \id) u^* \|  < \e,
\end{equation}
where $P_\s$ is the coordinate projection onto $\R^\s$,
and $M = C(\e) h \log h$.
We want to get a contradiction. 
To this end, note that the term $e_{h+1} \otimes e_{h+1}$ 
must be present in the expansion of $u P_\s u^*$ - 
in other words, $\nu$ must contain $N$.
Indeed, otherwise 
$$
\big( \a u P_\s u^* \big) e_{h+1}  =  0,
$$
although
$$
\| (u u^*) e_{h+1} \|  =  1.
$$
This will contradict (\ref{difference}).

Next put $\nu = \s \setminus \{ N \}$. 
Then, with $m_j = | A_j \cap \nu |$, \ $j = 1, \ldots, h$
we can write
\begin{eqnarray*}
u P_\s u ^* 
  & = &  \Big( \sum_{i \in \nu} x_i \otimes x_i \Big)  
         +  e_{h+1} \otimes e_{h+1}   \\  
  & = &  \Big( \sum_{j \le h}  \frac{m_j}{k} e_j \otimes e_j \Big)  
         +  e_{h+1} \otimes e_{h+1}.
\end{eqnarray*} 
Then by (\ref{ustaru}) we have
\begin{multline}                                \label{wanted}
\a u P_\s u^*  -  u u^* \\ 
   =   \Big[ \sum_{j \le h}  \Big( \a\frac{m_j}{k} -1 \Big)
                e_j \otimes e_j 
       \Big] 
         +  ( \a - 1) e_{h+1} \otimes e_{h+1}.
\end{multline} 
This is a diagonal operator, so its norm equals
$$
\max_{j = 1, \ldots, h} 
\Big\{  
  \Big| \a \frac{m_j}{k} - 1 \Big|, | \a - 1 |
\Big\}.
$$
Since by (\ref{difference}) this norm must be less than $\e$,
we have $\a < 1 + \e$.
On the other hand, 
$$
\sum_{j \le h} m_j  =  |\nu|  \le  n.
$$
Therefore there exists a $j \le h$ so that 
$m_j  \le  n / h$. Then 
$$
\Big| \a\frac{m_j}{k} - 1 \Big|
  \ge  1 - \Big| \a \frac{m_j}{k} \Big|
  \ge  1 - (1 + \e) \a \frac{n}{h k}
   =   1 - (1 + \e) \frac{n}{N}.
$$
If $0 < \e < 1/2$ and $h$ is essentially smaller than $N$, 
then $n / N$ is close to $0$, making the norm of 
operator in  (\ref{wanted}) close to $1$.
This contradicts (\ref{difference}).

\subsection{Tail probabilities.}

As a natural strengthening of Theorem \ref{marks},
we will compute the tail probability
$$                           
\Prob \{ \| u (\D - \id) u^* \| > t  \}.
$$
Since the operator $u u^*$ which we are approximating in 
Theorem \ref{marks} has norm at most one, the interesting 
range for $t$ is 
$0  <  t  <  1$.

\begin{proposition}                      \label{p:tail}
  For the random diagonal operator $\D$ in Theorem \ref{marks}, 
  letting 
  $$
  \e  =  C_0 \sqrt{\log n} \cdot \sqrt{ \frac{h}{n} },
  $$
  we have
  \begin{equation}                       \label{psitwo}
    \Prob \{ \| u (\D - \id) u^* \| > t  \}   
    \le  3 \exp \Big( - \frac{t^2}{\e^2} \Big)
  \end{equation}
  for all $0  <  t  <  1$.
\end{proposition}

This estimate is a consequence of the following 
lemma (\cite{Le-Ta} Lemma 3.7), 
which can be easily proved via expansion of the 
exponential function.
Recall the definition of the $\psi_p$-norm
of a random variable $Z$, for $p > 1$:
$$
\|Z\|_{\psi_p} 
  =  \inf \big\{ \l > 0 : \; 
                 \E \exp ( Z / \l )^p  \le  e
          \big\}.
$$  

\begin{lemma}          \label{easy}
  Let $Z$ be a positive random variable, and $d$ be an integer.
  Then the following are equivalent: 
  
  (i) $\|Z\|_p  \le  C_1 p^{d/2}$ for all $p \ge 1$.
  
  (ii)$\|Z\|_{\psi_{2/d}}  \le  C_2$.
  
  \noindent The constants $C_1, C_2$ only depend on each other.  
\end{lemma}
 
This observation reduces our problem to 
computing the $p$-th moment of $\| u (\D - \id) u^* \|$.
To this end note that the proof of Lemma \ref{Pisier}
gives an estimate on the $p$-th moment of 
$\big\| \sum_j  \e_j x_j \otimes y_j \Big\|$.
In particular case, for $x_j = y_j$, we have

\begin{lemma} (M.~Rudelson) \ \                            \label{Mark}
  Let $(y_j)$ be a finite set of vectors in $\R^m$. 
  Then for $p > 1$ 
  $$
  \Big( \E \Big\| \sum_j  \e_j y_j \otimes y_j  
     \Big\|^p 
  \Big)^{1/p}
  \le  C (p \vee \log m)^{1/2} \max_j \|y_j\|
       \cdot \Big\| \sum_j  y_j \otimes y_j \Big\|^{1/2}.
  $$
\end{lemma}

\begin{lemma}           \label{pnorms}
  Let $Z$ is a positive random variable 
  with $\|Z\|_\pone \ge 1$.
  Then 
  $$
  \|Z\|_p  \le  C p \log ( \E \exp Z )
  $$
  for all $p \ge 1$.
\end{lemma}

\proof
Let $M = \|Z\|_\pone \ge 1$, then 
$$
\E \exp (Z / M) = e.
$$
By Lemma \ref{easy}, $\|Z / M\|_p  \le  C p$.
Then by Jensen's inequality
\begin{eqnarray*}
\| Z \|_p
  &\le&  C p M 
    =    C p M \log ( \E \exp (Z /M) ) \\
  & = &  C p \log \Big( \E \exp (Z /M) \Big)^M \\
  &\le&  C p \log ( \E \exp Z).
\end{eqnarray*}
\endproof

\noindent {\bf Proof of Proposition \ref{p:tail}.\ }
Let 
$$
Z  =  u (\D - \id) u^*.
$$
Then \eqref{psitwo} is equivalent to 
\begin{equation}                          \label{zwedge}
  \big( \|Z\| \wedge 1 \big)_\ptwo
  \le  \e.
\end{equation}
By Lemma \ref{easy} there exists a constant $c > 0$ such that 
\eqref{zwedge} is implied by
$$
\big( \|Z\| \wedge 1 \big)_p  \le  c \e p^{1/2} 
\ \ \ \ \text{for all $p > 1$}.
$$
Note that 
$(\|Z\| \wedge 1)_p
 \le  \|Z\|_p \wedge 1$.
Then, with $E_p = \|Z\|_p$, it suffices to show that
\begin{equation}                          \label{pwedge}
  E_p \wedge 1  \le  c \e p^{1/2}
\end{equation}
for all $p > 1$.
Similarly to the proof of Theorem \ref{marks},
$$
E_p  \le  2 \Big( \E \E_\e 
               \Big\| \sum_{j=1}^N  
                      \e_j \D(j) x_j \otimes x_j  
               \Big\|^p
          \Big)^{1/p}.
$$
Applying Lemma \ref{Mark} in the same context as before
we get
\begin{eqnarray}                         \label{twoexp}
E_p  
&\le& C \left( 
          \E \left[ 
               \Big( p \vee \log \Rank \D \Big)^{p/2}
               \cdot K^{-p/2}
               \cdot \Big\| \sum_{j \le N} \D(j) x_j \otimes x_j 
                     \Big\|^{p/2}
             \right] 
        \right)^{1/p}        \nonumber\\
&\le& C K^{-1/2} 
        \Big[ \E \Big( p \vee \log \Rank \D \Big)^p
        \Big]^{1/2p}
        \Big[ \E \Big\| \sum_{j \le N} \D(j) x_j \otimes x_j 
                     \Big\|^p
        \Big]^{1/2p}    
\end{eqnarray}
To compute the first expectation we use Lemma \ref{pnorms}
for $Z = \log \Rank \D$.
Since $\Rank \D \ge e$ by the definition, $\|Z\|_\pone \ge 1$.
Therefore
\begin{eqnarray*}
\Big[ \E \Big( p \vee \log \Rank \D \Big)^p
\Big]^{1/2p}
  &\le&  ( p + \|Z\|_p)^{1/2}  \\
  &\le&  C ( p + p \log \E \Rank \D )^{1/2} \\
  &\le&  C \Big( p \log (K h) \Big)^{1/2} 
\end{eqnarray*}
as in the proof of Theorem \ref{marks}.
As for the second expectation in (\ref{twoexp}),
$$ 
\Big[ \E \Big\| \sum_{j \le N} \D(j) x_j \otimes x_j 
         \Big\|^p
\Big]^{1/2p} 
\le  (E_p + 1)^{1/2}.
$$
Therefore, recalling that $K = n / h$ we obtain
\begin{eqnarray}                  \label{david}
E_p  
  &\le&  C_0 \sqrt{\log n} \sqrt{ \frac{h}{n} } p^{1/2}
       (E_p + 1)^{1/2}   \notag\\
  &\le&  (c/10) \e p^{1/2} \cdot (E_p + 1)^{1/2},
\end{eqnarray}
provided the constant $C_0$ in the definition of $\e$
is large enough.

Now, if $(c/10) \e p^{1/2}  >  1$ then certainly
\eqref{pwedge} is true.
So we can assume that $(c/10) \e p^{1/2}  \le  1$.
Then \eqref{david} yields 
$$
E_p  \le  2  (c/10) \e p^{1/2}
     \le  c \e p^{1/2},
$$
which again implies \eqref{pwedge}. 
This proves \eqref{zwedge} and therefore completes the proof.
\endproof

\section{Isomorphisms on Coordinate Subspaces}

\subsection{Result}                                      \label{result}

In this section we will discuss the extension of 
Bourgain-Tzafriri's principle of restricted invertibility,
Theorem \ref{i:cor} and Corollary \ref{i:rip}, and 
its relation to the problem of harmonic density.

For the proof of Theorem \ref{i:cor} we refer to \cite{V}.
Note that it implies Corollary \ref{i:rip} by homogeneity.

One pleasant thing about Corollary \ref{i:rip} is that 
with {\em some} fixed $\e$ it can be deduced  
directly from Kashin-Tzafriri's suppression estimate \eqref{i:ktz} 
and the original Bourgain-Tzafriri's theorem \cite{B-Tz}, 
which we recall now.

\begin{theorem} (J.~Bourgain, L.~Tzafriri)              \label{btz}
  Let $u$ be a linear operator in $l_2^N$ such that 
  $\|u e_j\| = 1$ for all $j = 1, \ldots, N$.
  Then there exists a subset $\s$ 
  of $\{ 1, \ldots, N \}$ of cardinal 
  $|\s| > c \frac{N}{\|u\|^2}$
  so that
  $$                       
  \|u x\|  \ge  c \|x\|
    \ \ \ \ \text{for $x \in \R^\s$}.
  $$
\end{theorem}

To deduce Corollary \ref{i:cor}, first apply Kashin-Tzafriri's
suppression \eqref{i:ktz} to the operator $u$.
Note that $u$ is not necessarily norm one, so by homogeneity
we obtain a subset $\s_1$ of $\{1, \ldots, N \}$ of size
$|\s_1|  =  N / \|u\|^2$ so that the norm of $u$ on
$\R^{\s_1}$ is bounded by an absolute constant. 
Next apply Theorem \ref{btz} to $u$ restricted to $\R^{\s_1}$.
There exists a subset $\s$ of $\s_1$ 
of cardinality proportional to $|\s_1|$
so that $\|u x\|  \ge  c \|x\|$ for all $x \in \R^\s$.
Hence Corollary \ref{i:rip} is proved (for {\em some} $0 < \e < 1$).

The importance of the upper bound in Corollary \ref{i:rip} 
can be best illustrated by the following example,
which links this theme to harmonic analysis. 
Let $\T$ be the unit circle with the normalized Lebesgue 
measure $\nu$.
Consider an arc $B  \subset  \T$,
$$
B = \{ e^{it}, \ 0 \le t \le 2 \pi b \},
$$
and assume for simplicity that $\nu(B) = b$ is the inverse 
of a positive integer. 
Let $P_B$ be the restriction onto $B$:
$$
P_B f  =  f \chi_B.
$$
How does $P_B$ act on the natural "coordinate" structure 
generated by the characters $e^{i k t}$, $k \in \Z$?
$P$ maps the characters to vectors of norm $\sqrt{b}$
in $L_2 (\T)$.
Therefore, if $b$ is small then $P_B$ is far from being an
isomorphism on any "coordinate" subspace generated by 
the characters. 
On the other hand, an easy integration shows that
the vectors 
$$
P_B (e^{i k t}), \ \ k \in \frac{1}{b} \Z,
$$
are orthogonal.
Therefore, on the subspace of $L_2(\T)$ generated by 
the characters
$$
\Big\{ e^{i k t}, \ \ k \in \frac{1}{b} \Z \Big\}
$$
the projection $P_B$ is a multiple of an identity.

The finite-dimensional counterpart of this situation is
captured by Corollary \ref{i:rip}. Applied to the operator
$u = \frac{1}{\sqrt{b}} P_B$, it guarantees that the
set of the first $n$ characters contains 
a subset of size almost $b n$, 
on which $P_B$ acts like an isomorphism. 
The infinite-dimensional extension of this result is the 
subject of the next section.

\subsection{Harmonic Density Problem}     \label{HDP}

Our aim is to prove Theorem \ref{harmonic}.
Following \cite{B-Tz}, the proof will consist two steps.
First, the invertibility result, Corollary \ref{i:rip},
implies a finite-dimensional version of Theorem \ref{harmonic}.
Next we apply a combinatorial result of I.~Z.~Ruzsa \cite{Rus},
which allows to pass from large finite sets to a
needed infinite set $\L$ of large density.

\begin{definition}
  Let $\H$ be a set of finite sets of integers. 
  $\H$ is called a {\em homogeneous system} 
  if for every $A \in \H$, all the subsets 
  and translations of $A$ belong to $\H$. 
\end{definition}

\noindent Given a homogeneous system $\H$, there exists a limit
\begin{equation}                           \label{defdh}
d(\H)  =  \lim_{n \to \infty} \max_{A \in \H}
          \frac{ |A \cap [1, n]| } {n}.
\end{equation}
(note that the sequence under the limit is non-increasing).

\begin{theorem} (I.~Z.~Ruzsa \cite{Rus}). \           \label{ruzsa}
  Given an arbitrary homogeneous system $\H$,
  there exists a sequence of integers $\L$ such that
  its finite subsets all belong to $\H$ and
  $$
  \dens \L = d(\H).
  $$
\end{theorem}

In \cite{Rus} this theorem was proved for the  
one-sided density of $\L$, i.e. for 
$\lim_{n \to \infty} \frac{ |\L \cap [1, n]| } {n}$.
However in our case it seems more natural to work
with the two-sided density 
$$
\dens \L  =  \lim_{n \to \infty} 
             \frac{ |\L \cap [-n, n]| } {2n}.
$$
Our first buisness will be to prove Ruzsa's Theorem \ref{ruzsa}
for the two-sided density.
This requires only a slight modification of the original argument. 
At the first step, we find arbitrarily large finite 
subsets in $\H$ with optimal "hereditary density". 
This is the content of the following lemma 
from \cite{Rus}, which we include without proof. 

\begin{lemma} (I.~Z.~Ruzsa). \         
  Let $\H$ be a homogeneous system. 
  For arbitrary $n$ there exists a set 
  $A \in \H$,  $A \subset [1, n]$ satisfying
  $$
  \frac{| A \cap [1, k] |}{k}  \ge  d(\H)
  \ \ \ \ \text{for all $1 \le k \le n$}.
  $$
\end{lemma}

\noindent {\bf Proof of Theorem \ref{ruzsa}. \ }
By the homogeneity of $\H$, there also exists 
for arbitrary $n$ a set $A \in H$, $A \subset [-n, n]$ 
satisfying
\begin{equation}                   \label{hereditary}
  \frac{| A \cap [-k, k] |}{2 k}  \ge  d(\H)
  \ \ \ \ \text{for all $1 \le k \le n$}.
\end{equation}
We define an ordered tree $G$. Its $n$'th level $G_n$
consists of all the sets $A \in H$, $A \subset [-n, n]$ 
satisfying (\ref{hereditary}).
A vertice $A \in G_n$ is connected to a 
vertice $B \in G_{n+1}$ if $A \subset B$ .
This describes all the edges of $G$. 

Then the graph $G$ is infinite, 
joined (since every $A \in G_{n+1}$ is connected
to $A \cap [-n, n]  \in  G_n$), and every vertex of $G$
has only finitely many edges going to it. 
Then by K\"{o}nig's Infinity Lemma (see e.g. \cite{F})
there exists an infinite path in $G$.
So, let $(A_n)_{n \ge 1}$ be a chain in $G$:
$$
A_n \in G_n,  \ \ A_n \subset A_{n+1}, \ \ n = 1, 2, \ldots
$$
Then we claim that the set 
$$
\L = \cup A_n
$$
satisfies the conclusions of Theorem \ref{ruzsa}.
Indeed, all finite subsets of $\L$ belong to $\H$ by 
the construction.
To show that the density of $\L$ is $d(\H)$,
note that
$$ 
\frac{ |\L \cap [-n, n]| }{2 n}
\ge  \frac{ |A_n| }{2 n}
\ge  d(\H)
$$
from which it follows that the lower two-sided density 
$\underline{\dens} \L  \ge  d(\H)$.
Similarly, since the sets $\L \cap[-n, n]$ belong $\H$
and by homogeneity
$$
\frac{ |\L \cap [-n, n]| }{2 n}
\le  \max_{A \in \H} \frac{ |A \cap [-n, n]| }{2 n}
=    \max_{A \in \H} \frac{ |A \cap [1, 2n+1]| }{2 n}.
$$
Passing here to the upper limit as $n \to \infty$,
we conclude that the upper two-sided density 
$\overline{\dens} \L  \le  d(\H)$.
(note that in (\ref{defdh}) the limit exists).
This completes the proof of Theorem \ref{ruzsa}.
\endproof

\noindent {\bf Proof of Theorem \ref{harmonic}. \ }
Define an operator $u$ on $L_2(\T)$ as
$$
u f  =  \nu(B)^{-1/2} f \chi_B.
$$
Note that $\|u\|  =  \nu(B)^{-1/2}$ and
$\| u (e^{i k t}) \|  =  1$  for all $k \in \Z$.
Then we apply Corollary \ref{i:rip}.
For every positive integer $n$, we get a subset 
$$
\s_n  \subset  \{ 1, \ldots, n \}
$$
of cardinality 
$$
|\s_n|  \ge  (1 - \e) \nu(B) n
$$
for which
$$ 
c_1(\e) \|f\|_{L_2(\T)}  
\le  \nu(B)^{-1/2} \|f\|_{L_2(\T)}
\le  c_2(\e) \|f\|_{L_2(\T)}
$$
whenever the Fourier transform of $f$ is supported by $\s_n$.
The middle part of this inequality is $\|f\|_{L_2(B)}$.
Therefore we get:
$$
\text{Equivalence (\ref{equivalence}) holds 
      whenever $\supp \hat{f} \subset \s_n$.}
$$
Consider the family $\H$ of all finite subsets $\s$ 
of the integers such that the equivalence (\ref{equivalence}) 
holds whenever the Fourier transform of $f$ 
is supported by $\s$.
In particular, all sets $\s_n$ belong to $\H$.
Clearly, the family $\H$ is homogeneous.
Since $\s_n  \in  \H$, 
$$
d(\H)  
  \ge  \limsup_{n \to \infty} \frac{|\s_n|}{n}
  \ge  (1 - \e) \nu(B).
$$
Then Ruzsa's Theorem \ref{ruzsa} yields the existence of a
set $\L$ of integers whose all finite subsets
belong to $\H$, and with two-sided density
$$
\dens \L  \ge  (1 - \e) \nu(B).
$$
This completes the proof in view of the definition of $\H$.
\endproof

To see how sharp Theorem \ref{harmonic} is, let us look again 
at the example in Section \ref{result}.
We consider an arc
$$
B = \{ e^{it}, \ 0 \le t \le 2 \pi b \},
$$
as a subset of $\T$, where $\nu(B) = b$ is the inverse 
of a positive integer. 
We will show that if a set of integers $\L$ has two-sided
density exceeding $\nu(B)$, then there exists a function $f$
with $\supp \hat{f}  \subset \L$ and for which 
\eqref{equivalence} fails; more precisely
\begin{equation}                                 \label{fails}
  \|f\|_{L_2(\T)}  =  1
  \ \ \ \ \text{and} \ \ \ \ 
  \|f\|_{L_2(B)}  \le \a
\end{equation}
where $\a$ can be chosen arbitrarily small.

Tu put this differently, \eqref{fails} means that 
the sequence $\{ e^{i k t}, k \in \L \}$
is not a Riesz basis in $L_2(B)$. 
Assume the opposite. Then $\{ e^{i k t}, k \in \L \}$
is equivalent to the canonical basis in $l_2$, i.e.
$$
\int_0^b  \Big| \sum_k a_k e^{2 \pi i k t} \Big|^2 \; dt
\sim \sum_k |a_j|^2
$$
for all finite sets of scalars $(a_k)$.
By change of variable,
$$
\int_0^1  \Big| \sum_k a_k e^{2 \pi i b k t} \Big|^2 \; dt
\sim \sum_k |a_j|^2,
$$
showing that for $\{ \l_k \} := b \L$, the sequence
of exponentials $(e^{i \l_k t})$ is equivalent
in $L_2(\T)$ to the canonical basis of $l_2$.
Then we apply a classical result of N.~Levinson
on completeness of exponentials in $L_p(\T)$,
see \cite{Lev} Appendix III.1 or \cite{You} 3.2.
Let $n(r)$ denote the number of points $\l_k$ inside
the disc $|z|  \le  r$, and we put
$$
N(r)  =  \int_1^r  \frac{n(t)}{t} \; dt.
$$

\begin{theorem} (N.~Levinson). \                      \label{levinson}
  The set $(e^{i \l_k t})$ is complete in $L_p(\T)$ 
  whenever
  $$
  \limsup_{r \to \infty} 
    \Big( N(r) - 2 r + \frac{1}{p} \log r \Big)
  > \infty.
  $$
\end{theorem}

In our setting the ration $\frac{n(t)}{2 t}$ approaches
$\frac{1}{b} \dens \L  >  1$ as $N \to \infty$.
Thus we have $N(r) > 2 r$ for $r$ sufficiently large. 
Then by Theorem \ref{levinson} the system $(e^{i \l_k t})$
is complete in $L_2(\T)$.

Note that this argument holds also if we remove a finite 
number of elements from $(e^{i \l_k t})$,
so that this system remains complete after the removal.
This clearly contradicts to its equivalence to the 
canonical basis in $l_2$. This finishes the proof.

\remark The use of the result of N.~Levinson was suggested 
to me by V.~Kadets.

\section{APPENDIX. Application to Communication Systems}    

We apply results of Section \ref{global} to provide 
optimal estimates for a communication system which delivers
data with random losses \cite{G-K}.

A typical case we have in mind here is the Internet. 
A requested information is sent to a user in a sequence of 
"data packets". 
If a data packet is lost on its way to the user, a protocol detects 
the missing packet and sends it again. 
However, the detection of the lost packet usually takes much more 
time than a successful delivery.
This is the main source of large delays known to all network users.
Therefore, instead of retransmitting the lost packets
it is highly desirable to be able to recover the sent information 
using whatever received, despite the loss of some packets.
The question is then how to distribute the source information 
among data packets?
There should be some dependency between the packets,
otherwise the information contained in the missing packets is 
irrevocably lost.

Parallel to the development of wavelets and connected with it,
there has arisen a simple but fruitful idea to represent 
information, viewed as a vector $x$ in $\R^m$, by its expansion 
through an identity 
\begin{equation}                                     \label{frame}
  id = \sum_{j \le k} x_j \otimes x_j
\end{equation} 
for suitable vectors $x_j  \in  \R^m$. 
These vectors are called a {\em frame} \cite {Da}. 
(More generally, a frame is a set of vectors $x_j$ for which 
$\sum x_j \otimes x_j$ is an isomorphism in $l_2^m$).
Clearly, $k \ge m$.

This way, a source vector $x$ in $\R^m$ is represented 
by $k$ data packets -- coefficients 
$\< x_j, x\> $, $j = 1, \ldots, k$, which carry complete 
information about $x$ due to the reconstruction formula
\begin{equation}                       \label{reconstruction}
  x  =  \sum_{j\le k} \< x_j, x\> x_j.
\end{equation}
If $k > m$ then this information is redundant; there is a 
kind of dependency between the packets.
This way of representing $x$ is often more resiliant 
to errors than the old method -- expanding $x$ using
an orthonormal basis and transmitting each coefficient 
$k / m$ times \cite{Da}.

A problem raised in \cite{G-K} was: is this new method 
also resiliant to random losses of the packages?
Specifically, if a random (but not too large) subset 
of the packages $\< x_j, x\> $, $j = 1, \ldots, k$
is lost on the way to a user, can one essentially recover $x$ 
by summing in \eqref{reconstruction} only the successfully
delivered components?

To make this scheme work with probability at least $1/2$, 
the norms $\|x_j\|$ have to be reasonably small, otherwise
the contribution of the summands in \eqref{reconstruction}
can be too irregular; one can easily produce examples 
making this intuitive statement precise.

Under this assumption the results of Section \ref{global} 
imply that if at least $C(\e) m \log m$ packets
are successfully delivered, then with high probability 
the source vector $x$ can be 
reconstructed by \eqref{reconstruction} with precision $\e$. 
The two main points here are that

\begin{itemize}
  \item{The required number of successfully delivered packets
   does not depend on $k$;}
  \item{no information is needed about the lost packets.}
\end{itemize} 

As for the mentioned restriction on the norms $\|x_j\|$,
we will assume for simplicity they are equal to each other 
(and therefore to $\sqrt{ \frac{m}{k} }$);
a more general case requires only minor changes.

\begin{theorem}
  Consider a set of vectors $x_1, \ldots, x_k$ in $\R^m$ 
  with equal norms and which satisfy \eqref{frame}. 
  Let $\s$ be a random subset of $\{1, \ldots, k\}$ 
  with cardinal $|\s|  \sim  n$ (i.e. each element
  of $\{1, \ldots, k\}$ is taken or rejected independently 
  with probability $n/k$).
  Then for $0 < t < 1$
  \begin{equation}                          \label{lesst}
    \Big\| \id  -  \frac{k}{|\s|} 
                   \sum_{j \in \s} x_j \otimes x_j
    \Big\|  <  t
  \end{equation}
  with probability at least $1 - 6 \exp \big( - t^2 / \e^2 \big)$,
  where $\e  = C \sqrt{\log n} \cdot \sqrt{ \frac{m}{n} }$.
\end{theorem}

\proof
Note that $\|x_j\| = \sqrt{ \frac{m}{k} }$ for all $j$.
We will apply Theorem \ref{marks} to the operator 
$u: l_2^k \to l_2^m$ defined by 
$$
u e_j  =  x_j,    \ \ \ \ j = 1, \ldots, k.
$$
It can easily be seen that 
$$
\|u\| = 1
\ \ \ \ \text{and} \ \ \ \
h = \|u\|_\HS^2 = m.
$$
Also
$$
u u^*  =  \sum_j x_j \otimes x_j  =  \id, \ \ \ \
u \D u^*  =  \sum_j \D(j) x_j \otimes x_j,
$$
where $\D(j)$ denotes the $j$-th diagonal entry of $\D$.
Recall that, by the construction of $\D$, in our case
$\D(j)$ is a random variable independent of the other 
entries of $\D$ and distributed as
$$
\Prob \Big\{ \D(j) = \frac{k}{n} \Big\}  
  =  1 - \Prob \big\{ \D(j) = 0 \big\}
  =  \frac{n}{k}.
$$ 
Then 
\begin{equation}                    \label{adone}
\| u (\D - \id) u^* \|  
  =  \Big\| \id  -  \frac{k}{n} 
                   \sum_{j \in \s} x_j \otimes x_j
     \Big\|   
\end{equation}
where $\s = \{  j: \ j \le k, \D(j) \ne 0 \}$
is a random subset of $\{ 1, \ldots, N \}$,
as described in the assumption of the theorem.

Then Theorem \ref{marks} and Proposition \ref{p:tail}
give the needed probability of \eqref{lesst},
but with $|\s|$ replaced by $n$ in the denominator.
(However, we do not want $n$ in the reconstruction formula
since a priori $n$ is not known). 
To complete the proof we note that 
by a standard concentration inequality
$|\s|$ is close to $n$ with high probability.
Indeed, $|\s|$ is a sum of $k$ independent $\{ 0, 1\}$-valued
random variables $\d_j$ with $\d := \E \d_j = \frac{n}{k}$.
Then it follows from the classical bounds on the binomial law
\cite{Ho} that for $s  \le  2 \d k$
$$
\Prob \big\{  \big| |\s| - n \big| > s \big\}
 \le  \exp \Big( \frac{s^2}{4 \d (1 - \d) k}
            \Big)   
 \le  \exp \Big( - \frac{s^2}{8 n} \Big).
$$
Then for $s = t n$
$$
\Prob \Big\{  \Big| \frac{|\s|}{n} - 1 \Big| > t \Big\}
 \le  \exp \Big( - \frac{t^2 n}{8} \Big)    
 \le  \exp \Big( - \frac{t^2}{\e^2} \Big).     
$$
This justifies the replacement of $|\s|$ by $n$ 
and therefore completes the proof.
\endproof

\end{document}